\documentclass[reqno, 12pt]{amsart}

\usepackage{amssymb}
\usepackage{latexsym}
\usepackage{epsfig}
\usepackage{epsf}
\usepackage{float}
\usepackage{a4wide}
\usepackage{hyperref}

\newtheorem{theorem}{Theorem}
\newtheorem{lemma}[theorem]{Lemma}
\newtheorem{corollary}[theorem]{Corollary}

\def\bR{{\mathbb R}}  
\def\bN{{\mathbb N}}  
\def\bP{{\mathbb P}}    
\def\bZ{{\mathbb Z}}

\def\bE{{\mathbb E}}

\def\N{{\mathcal N}}

\def\V{{\mathcal V}}

\def\calC{{\mathcal C}}

\def\calN{{\mathcal N}}
\def\calX{{\mathcal X}}
\def\calY{{\mathcal Y}}
\def\calZ{{\mathcal Z}}
\def\calS{{\mathcal S}}

\def\D{{\mathcal D}}

\def\P{{\mathcal P}}

\def\A{{\mathbf A}}
\def\B{{\mathbf B}}
\def\K{{\mathbf K}}
\def\C{{\mathbf{Cl}}}

\def\K{{\mathbf K}}
\def\0{{\mathbf 0}}
\def\1{{\mathbf 1}}
\def\e{{\mathbf e}}

\def\v{{\mathbf v}}
\def\u{{\mathbf u}}
\def\x{{\mathbf x}}
\def\y{{\mathbf y}}
\def\z{{\mathbf z}}

\def\II{\mbox{ 1\hskip -.29em I}}

\def\reff#1{(\ref{#1})}
\def\proofof #1{{\noindent \emph{Proof of #1}.}}
\def\endproof{{\flushright $\square$}} 

\begin{document}

\title[On some fundamental aspects of Polyominoes on Random Voronoi Tilings]{On Some fundamental aspects of Polyominoes\\ on Random Voronoi Tilings}
\author{Leandro P. R. Pimentel}
\address{Institute of Mathematics\\
Federal University of Rio de Janeiro\\ 
Caixa Postal 68530, CEP 21941-909, Rio de Janeiro, RJ, Brazil.} 
\email{leandro@im.ufrj.br}
\urladdr{}



\begin{abstract}
Consider a Voronoi tiling of $\bR^d$ based on a realization of a inhomogeneous Poisson random  set. A Voronoi polyomino is a finite and connected union of Voronoi tiles. In this paper we provide tail bounds for the number of boxes that are intersected by a Voronoi polyomino, and vice-versa. These results will be crucial to analyze self-avoiding paths, greedy polyominoes and first-passage percolation models on Voronoi Tilings and on the dual graph, named the Delaunay triangulation \cite{P07,PR08}. 
\end{abstract}

\maketitle

\section{Introduction}
To any locally finite subset $\N$ of $\bR^d$ one can associate a
partition of the plane as follows. To each point $\v\in \N$
corresponds a polygonal region $C_\v$, the \emph{Voronoi tile} at $\v$, consisting of the set of points of $\bR^d$ which are
closer to $\v$ than to any other $\v'\in \N$. Closer is understood here in the euclidean sense, and
the partition is not a real one, but the set of points which belong to
more than one Voronoi tile has Lebesgue measure 0. From now on, $\N$ is
understood to be distributed like a Poisson random set on $\bR^d$ with
intensity measure $\mu$. We shall always assume that $\mu$ is
comparable to Lebesgue's measure on $\bR^d$, $\lambda_d$, in the sense
that there exists a positive constant $c_\mu$ such that for every
Lebesgue-measurable subset $A$ of $\bR^d$:
\begin{equation}
\label{comparaison}
c_\mu^{-1}\lambda_d(A)\leq\mu(A)\leq c_\mu\lambda_d(A)\;.
\end{equation}
Notice that, with probability one,
when two Voronoi tiles are connected, they share a $(d-1)$-
dimensional face. The collection
$\V=\V(\N):=\{C_\v\}_{\v\in \N}$ is called the \emph{Voronoi tiling} (or
tessellation) of the plane based on $\N$. 

\begin{figure}[htb]
\begin{center}
\includegraphics[width=0.6\textwidth]{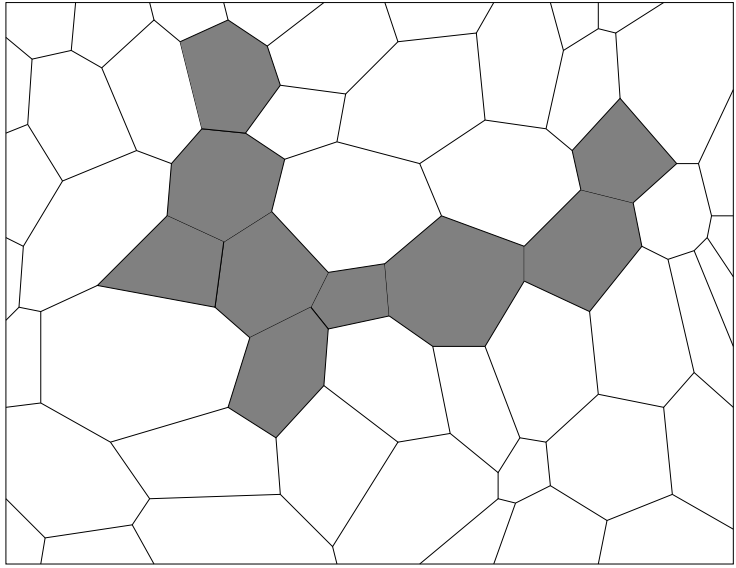}
\caption{A two-dimensional Voronoi tiling and a (gray colored) Voronoi polyomino of size $n=9$ .}\label{poly}
\end{center}
\end{figure}

The study of Voronoi tilings has a very long history. The terminology is in honour of Voronoi \cite{V08}, who used these tilings to study quadratic forms. Our aim is to study some fundamental aspects of finite and connected union of Voronoi tiles. These objects are called \emph{Voronoi polyominoes} (Figure \ref{poly}). Polyominoes were first introduced in periodic tilings of the plane. The Voronoi setting is rather different from the periodic one since we are now considering a random environment induced by the underlying Poisson random set. 

The main results of this article, that will be stated in Section \ref{Po}, will provide tail bounds for the maximum and minimum number of square boxes intersected by a Voronoi polyomino that contains the origin $\0$ and has size $r$. They will be important tools to analyze self-avoiding paths, greedy polyominoes and first-passage percolation models on the Voronoi random setting \cite{P07,PR08}.  The idea to prove them is to combine block arguments with standard results for greedy lattice animal and site percolation models. The block argument is to consider a large box in $\bR^d$ so that it contains with ``high probability'' some configuration of points which prevent a Voronoi tile to cross it completely. The ``high probability'' alluded to is some percolation probability: we need that the ``bad boxes'' (those who can be crossed) do not percolate. In Section \ref{Te} we will prove some technical lemmas concerning greedy lattice animals and site percolation models that will be combine with the block argument to prove the theorems in Section \ref{Pr}. In Section \ref{Mo} a slightly different random set up is introduced and we will draw an outline of how the results can be extended to this new set up.   
       
\section{Polyominoes on Voronoi Tilings}\label{Po}
Let $\#S$ denote the usual cardinality of a set $S$, and for each subset $A$ of $\bR^d$ let $\#_\N A:=\#(A\cap\N)$. For each natural number $r\geq 1$ a Voronoi \emph{
polyomino} $\P$ of size $r$ is a connected union of $r$
Voronoi tiles. Let $\Pi_{\geq r}$ denote the collection of all polyominoes $\P$ such that the origin $\0\in\P$ and $\#_\N\P\geq r$, and let $\Pi_{\leq r}$ be the collection of all polyominoes $\P$ such that $\0\in\P$ and $\#_\N\P\leq r$. Let $\bZ^d$ denote the $d$-dimensional integer lattice. For each $\z\in\bZ^d$ let
$$B_\z:=\z+[-1/2,1/2)^d\,,$$
and for each connected set $C\subseteq\bR^d$ let
$$\A(C):=\left\{\z\in\bZ^d\,\,:\,\,B_\z\cap C\neq\emptyset\,\right\}\,.$$

\begin{theorem}\label{t1}
There exists a constant $b_1\in(0,\infty)$ such that if $\,r\geq b_1 s\,$ then 
\begin{equation}\label{e1t1}
\bP\left( \min_{\P\in\Pi_{\geq r}}\#\A(\P)\leq s\right)\leq e^{-r/2}\,.
\end{equation}
Further, there exist constants $b_2,b_3\in(0,\infty)$ such that  if $\,s\geq b_2 r\,$  then 
\begin{equation}\label{e2t1}
\bP\left(\max_{\P\in\Pi_{\leq r}}\#\A(\P)\geq s\right)\leq e^{-b_3 s}\, .
\end{equation}
\end{theorem}

In the same polyomino model one could consider  $\max$ and $\min$ of $\#\A(\P)$ over all polyominoes $\P$ of size $r$ touching $B_\0$ ($\P\cap B_\0 \neq \emptyset$). The same method to prove Theorem \ref{t1} can be extended to this situation, yielding to similar large deviations bounds \eqref{e1t1} and \eqref{e2t1}.

\subsection{Self-avoiding paths on the Delaunay Triangulation}
\begin{figure}[tb]
\begin{center}
\includegraphics[width=0.6\textwidth]{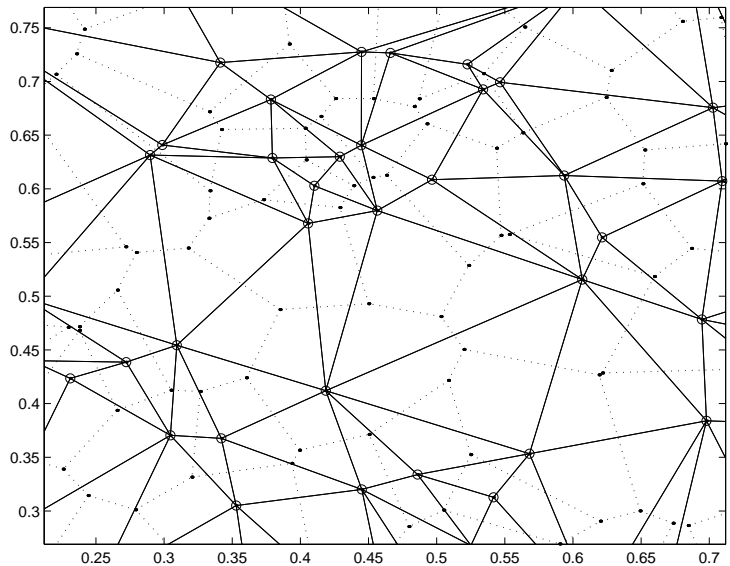}
\end{center}
\caption{ The Voronoi Tiling $\V$ (dotted lines) and the Delaunay Triangulation $\D$ (solid lines) in the two-dimensional model.}\label{vor}
\end{figure}

An important graph for the study of a Voronoi tiling is its facial dual, the \emph{Delaunay graph} based on $\N$. This graph, denoted by $\D=\D(\N)$ is an unoriented graph embedded in $\bR^d$ which has vertex set $\N$ and edges $\{\u,\v\}$ every time $C_\u$ and $C_\v$ share a $(d-1)$-dimensional face. We remark that, for our Poisson random set, a.s. no $d+1$ points are on the same hyperplane and no $d+2$ points are on the same hypersphere and that makes the Delaunay graph a well defined triangulation (Figure \ref{vor}). This triangulation divides $\bR^d$ into bounded simplexes called  \emph{Delaunay cells}. For each Delaunay cell $\Delta$ no point in $\N$ is inside the circum-hypersphere of $\Delta$. Polyominoes on the Voronoi tiling correspond to connected (in the graph topology) subsets of the Delaunay graph.

Let $\v_\x$ be the nearest point $\v\in\N$ to $\x\in\bR^d$, and let $\Gamma_{\geq r}$ (resp., $\Gamma_{\leq r}$) be the collection of all self-avoiding paths $\gamma$ starting at $\v_\0$ and of size (number of vertices) $\#\gamma\geq r$ (resp.,  $\#\gamma\leq r$). Recall that polyominoes on the Voronoi
tiling correspond to connected (in the graph topology) subsets of the Delaunay graph. Therefore, for each $\gamma\in\Gamma_{\geq r}$ corresponds a unique polyomino $\P_\gamma\in\Pi_{\geq r}$ and, analogously, for each $\gamma\in\Gamma_{\leq r}$ corresponds a unique polyomino $\P_\gamma\in\Pi_{\leq r}$. Let $\A(\gamma):=\A(\P_\gamma)$. Thus, the following corollary is a straightforward consequence of Theorem \ref{t1}.
\begin{corollary}\label{c1}  
There exists a constant $b_1\in(0,\infty)$ such that if $\,r\geq b_1 s\,$ then
\begin{equation}\label{e1c1}
\bP\left( \min_{\gamma\in\Gamma_{\geq r}}\#\A(\gamma)\leq s\right)\leq e^{-r/2}\,.
\end{equation}
Further, there exist constants $b_2,b_3\in(0,\infty)$ such that if $\,s\geq b_2 r\,$  then
\begin{equation}\label{e2c1}
\bP\left(\max_{\gamma\in\Gamma_{\leq r}}\#\A(\gamma)\geq s\right)\leq e^{-b_3 s}\, .
\end{equation}
\end{corollary}

 \subsection{The inverse problem} Until now we have been concerned with the size of a lattice covering of a Voronoi polyomino of size $r$. It is also natural to consider the inverse problem, i.e., the number of Voronoi tiles that one needs to cover a connected set composed by $s$ lattice boxes. Precisely, for each connected set $\A\subseteq\bZ^d$ (in the $l_1$-nearest-neighbor sense) let 
$$B_\A:= \bigcup_{\z\in\A} B_\z\,.$$ 
 A connected subset of $\bZ^d$ is also called a \emph{lattice animal}. We denote $\Phi_{\leq s}$  the collection of all lattice animals such that  $\0\in\A$ and $\#\A\leq s$. A Voronoi covering is defined by taking 
 $$\P(\A):=\left\{\x\in\bR^d\,\,:\,\,\x\in C_\v\mbox{ and }C_\v\cap B_\A\neq\emptyset\,\right\}\,.$$ 
\begin{theorem}\label{t2}
There exist constants $b_7,b_8\in(0,\infty)$ such that if $\,r\geq b_7 s\,$  then
\begin{equation}\label{e2t2}
\bP\left(\max_{\A\in\Phi_{\leq s}}\#_\N\P(\A)\geq r\right)\leq 2e^{-b_8 r}\, .
\end{equation}
\end{theorem}

One important consequence of Theorem \ref{t2} concerns the following construction: Let $\x,\y\in\bR^d$ and consider the Polyomino $\P([\x,\y])$ generated by all Voronoi tiles that intersect the line segment $[\x,\y]$. Then one can always find a self-avoiding path $\gamma(\x,\y)$ with vertices in $\P([\x,\y])$ and that connects $\v_\x$ to $\v_\y$. Clearly,
$$\#\gamma(\x,\y)\leq \#_\calN\P([\x,\y])\,,$$
and hence, by Theorem \ref{t2}, we have the following corollary:
\begin{corollary}\label{c2}  
There exist constants $b_7,b_8\in(0,\infty)$ such that if $\,r\geq b_7 s\,$  then
\begin{equation}\label{e1c2}
\bP\left(\max_{\x:\|\x\|_2\leq s}\#\gamma(\0,\x)\geq r\right)\leq 2e^{-b_8 r}\,,
\end{equation}
where $\|\,.\,\|_2$ denotes the euclidean norm.
\end{corollary} 

\section{Technical Lemmas}\label{Te}
\subsection{A greedy lattice animal model with Poisson weights} 
Let 
$$N_\z= \#_\N B_\z\,.$$ 
Then $\{N_\z\,\,:\,\,\z\in\bZ^d\}$ is a collection of independent Poisson random variables. By \reff{comparaison},
\begin{equation}\label{exp-comparison}
\sup_{\z\in\bZ^d}\bE(e^{N_\z})\leq e^{c_\mu(e-1)}\,.
\end{equation} 
It is a standard result in combinatorics\footnote{To see this, notice that for each lattice animal $\A\in\Phi_{\leq s}$ one can (injectively) associate an ``exploration'' nearest neighbor path $(\0,\z_1,\dots,\z_l)$ such that $\z_i\in\A$ and $\#\{\z_i\,:\,\z_i=\z\}\leq 2d$ for each $\z\in\A$. Thus  $l\leq 2ds$ and $\alpha:=(2d)^{2d}$ will do.} that
\begin{equation}\label{counting}
\#\Phi_{\leq s}\leq \alpha^s\,,
\end{equation}
for a finite $\alpha=\alpha(d)$.
\begin{lemma}\label{Poisson}
If $r\geq 2(\log \alpha +c_\mu(e-1))s$ then
$$ \bP\left(\max_{\A\in\Phi_{\leq s}}\sum_{\z\in\A} N_\z\geq r\right)\leq e^{-r/2}\,.$$
\end{lemma} 
\proofof{Lemma \ref{Poisson}}
Combining \reff{exp-comparison} and \reff{counting} together with Markov's inequality, one has that
\begin{eqnarray}
\nonumber \bP\left(\max_{\A\in\Phi_{\leq s}}\sum_{\z\in\A} N_\z\geq r\right)&\leq&\sum_{\A\in\Phi_{\leq s}}\bP\left(\sum_{\z\in\A} N_\z\geq r\right)\\
\nonumber&\leq&\alpha^s\left[\sup_{\z\in\bZ^d}\bE(e^{N_\z})\right]^s e^{-r}\\
\nonumber&\leq& \exp\left\{\left[\log \alpha +c_\mu(e-1)\right]s-r\right\}\leq e^{-r/2}\,,
\end{eqnarray}
whenever $r\geq 2(\log \alpha +c_\mu(e-1))s$.
\begin{flushright}\endproof\end{flushright} 

\subsection{Site percolation schemes}
Throughout this section $\calY:=\{Y_{\z}\,\,:\,\,\z\in\bZ^d\}$ will denote an i.i.d. site percolation scheme (or random field) with parameter $\rho\in(0,1)$. If $Y_\z=1$ we say that $\z$ is open. Otherwise, we say that it is closed.  
 
\begin{lemma}\label{opendensity}
If $2\alpha\sqrt{1-\rho}<e^{-1}$ and $s\geq 2r$ then
$$ \bP\left(\min_{\A\in\Phi_{\geq s}}\sum_{\z\in\A} Y_\z \leq r\right)\leq e^{-s}\,,$$
where $\Phi_{\geq s}$ denotes the set of all connected sets $\A\subseteq\bZ^d$ such that $\#\A\geq s$ and $\0\in\A$.
\end{lemma}

\proofof{Lemma \ref{opendensity}} 
Let $\binom{s}{r}$ denote the binomial coefficient. If $\A\in\Phi_{\geq s}$ and $\sum_{\z\in\A} Y_\z \leq r $, by taking a connected subset of $\A$ of size $s$ we may assume that $\A$ has size exactly $s$. Then there exists some subset of exactly $s-r$ sites of $\A$ with $Y_\z=0$. By \reff{counting}, this shows that
\begin{eqnarray}
\nonumber\bP\left(\min_{\A\in\Phi_{\geq s}}\sum_{\z\in\A} Y_\z \leq r\right)&\leq&\alpha^s\binom{s}{r}(1-\rho)^{s-r}\\
\nonumber&\leq&\alpha^s2^s(1-\rho)^{s-r}\\
\label{gooddensity}&\leq& (2\alpha\sqrt{1-\rho})^{s}\leq e^{-s}\,,
 \end{eqnarray}
 whenever $2\alpha\sqrt{1-\rho}<e^{-1}$ and $s\geq 2r$.
 \begin{flushright}\endproof\end{flushright}   

 A closed cluster is a maximal connected set of closed vertices of $\bZ^d$. Let $\C_\z$ denote the closed cluster that contains $\z$ (it is empty if $\z$ is open). For a finite subset $\A$ of $\bZ^d$  let $\calC_\calY(\A)$ denote the collection of all closed clusters (with respect to $\calY$) intersecting $\A$. 
\begin{lemma}\label{lemcorrelneg}
 If $f$ is an increasing function from $\bN$ to $[1,+\infty[$ then
$$\bE\left(\prod_{\C\in\calC_\calY(\A)}f(\#\C)\right)\leq \left\{\bE f(\#\C_\0)\right\}^{\#\A}\,.$$
\end{lemma}
\proofof{Lemma \ref{lemcorrelneg}} The proof of Lemma \ref{lemcorrelneg} is due to Rapha\"{e}l Rossignol and the author is grateful for his help. This inequality is also a crucial tool for proving the results in \cite{PR08}. Let us recall Reimer's inequality (see Grimmett's book \cite{G99}
p.39). Let $n$ be a positive integer, let $\B(n)=\bZ^d\cap [-n,n]^d$
and define $\Omega_n =\{0,1\}^{\B(n)}$. For $\omega\in \Omega_n$ and
$\K\subset \B(n)$, define the cylinder event $C(\omega,\K)$ generated
by $\omega$ on $\K$ by:
$$C(\omega,\K)=\{\omega'\in\Omega_n\mbox{ s.t. }\omega'_i=\omega_i\forall i\in
\K\}\;.$$ If $A$ and $B$ are two subsets of $\Omega_n$, define their
disjoint intersection $A\square B$ as follows:
$$
A\square B=\{\omega\in\Omega_n\,:\,\exists\, \K\subset
\B(n),\;C(\omega,\K)\subset A\mbox{ and }C(\omega,\K^c)\subset
B\}\;.
$$
Reimer's inequality states that:
$$\bP(A\square B)\leq \bP(A)\bP(B)\;.$$
Remark that $\square$ is a commutative and associative operation,
and that, for any $l$ subsets $A_1,\ldots ,A_l$ of $\Omega_n$,
\begin{eqnarray*}
A_1\square \ldots \square A_l&=&\Big\{\omega\in\Omega_n\,:\,\exists\, \K_1,\ldots,\K_l\mbox{ disjoint subsets of
  }\B(n),\\
&&\;\bigcup_{i=1}^l\K_i=\B(n)\mbox{ and }C(\omega,\K_i)\subset
A_i\forall
  i=1,\ldots ,l\Big\}\;.
\end{eqnarray*}
Now take $n$ large enough so that $\A\subset \B(n)$. Let
$l=\#\A$, and order the elements of
$\A=\{\x_1,\x_2,\ldots,\x_l\}$. For any $\x$ in $\B(n)$, and
any $\omega\in\Omega_n$, let $\C_n(\x,\omega)$ be the closed
cluster (for the configuration $\omega)$ in $\B(n)$ containing $\x$,
which is empty if $\omega(\x)=1$. Define:
$$
\forall i,\;
\C_i(\omega)=\left\lbrace\begin{array}{ll}\C_n(\x_i,\omega)
&\mbox{ if }\x_i\not\in\bigcup_{k=1}^{i-1}\C_n(\x_k,\omega) \\
\emptyset &\mbox{
  else.}\end{array}\right.\;
$$
Let $k_1,\ldots ,k_l$ be nonnegative integers and:
$$
A_i=\{\omega\,:\,\#\C_i(\omega)\geq k_i\},\;\forall 1\leq
i\leq l\;,
$$
$$
\widetilde{A}_i=\{\omega\,:\,\#\C_n(\x_i,\omega)\geq
k_i\},\;\forall 1\leq i\leq l\;.
$$
Then, we claim that:
$$
\bigcap_{i=1}^l A_i\subset \widetilde{A}_1\square \ldots \square
\widetilde{A}_l\;.
$$
Indeed, let $\omega\in \bigcap_{i=1}^l A_i$. Then, for every
$i=1\ldots ,l-1$, $C(\omega,\C_i(\omega))\subset A_i$.
Furthermore, $C(\omega,\B(n)\setminus
\bigcup_{i=1}^{l-1}\C_i(\omega))\subset A_l$. This shows that
$\omega\in \widetilde{A}_1\square \ldots \square \widetilde{A}_l$,
and proves the claim above. Therefore, using Reimer's inequality,
\begin{equation}
\label{ineqReimer}\bP\left(\bigcap_{i=1}^l
A_i\right)\leq\Pi_{i=1}^l\bP\left(\widetilde{A}_i\right)\;.
\end{equation}
Now, let $f$ be an increasing function from $\bN$ to $[1,\infty)$.
Define $f_1$ as follows:
$$
f_1(0)=1\mbox{ and }\forall k\geq 1,\;f_1(k)=f(k)\;.
$$
Denote by $\{\beta_0=1<\beta_1 <\ldots <\beta_j<\ldots \}$ the
range of $f_1$. Define:
$$
k_j=\inf\{k\,:\,f(k)\geq \beta_j\},\;\forall j\in\bN\;,
$$
$$
A_{i,j}=\{\omega\,:\,\#\C_i(\omega)\geq
k_j\}=\{\omega\,:\,f_1(\#\C_i(\omega))\geq \beta_j\}\;,
$$
and
$$
\widetilde{A}_{i,j}=\{\omega\,:\,\#\C_n(\x_i,\omega)\geq
k_j\}=\{\omega\,:\,f_1(\#\C_n(\x_i,\omega))\geq
\beta_j\}\;.
$$
By convention, set $\beta_{-1}=0$ and define
$a_j=\beta_j-\beta_{j-1}$. We can write:
$$
f_1(\#\C_i(\omega))=\sum_{j\in\bN}(\beta_j-\beta_{j-1})\II_{A_{i,j}}=\sum_{j\in\bN}a_j\II_{A_{i,j}}\;.
$$
Define $\calC_n(\A)$ as the set of nonempty components
$\C_n(\x_i,\omega)$, for $i\in\{1,\ldots ,l\}$. Since $f_1(0)=1$,
we can write:
\begin{eqnarray*}
\bE\left(\Pi_{\C\in\calC_n(\A)}f(\#\C)\right)&=&\bE\left(\Pi_{\C\in\calC_n(\A)}f_1(\#\C)\right)\\
&=&\bE\left(\Pi_{i=1}^lf_1(\#\C_i(\omega))\right)\\
&=&\bE\left(\Pi_{i=1}^l\sum_{j\in\bN}a_j\II_{A_{i,j}}\right)\\
&=&\sum_{j_1,\ldots ,j_l}a_{j_1}\ldots
  a_{j_l}\bE\left(\Pi_{i=1}^l\II_{A_{i,j_i}}\right)\\
&\leq &\sum_{j_1,\ldots ,j_l}a_{j_1}\ldots
  a_{j_l}\Pi_{i=1}^l\bP\left(\widetilde{A}_{i,j_i}\right)\\
&=&\Pi_{i=1}^l\sum_{j\in\bN}a_j\bP\left(\widetilde{A}_{i,j}\right)\\
&=&\Pi_{i=1}^l\bE\left(f_1(\#\C_n(\x_i,\omega))\right)\leq\Pi_{i=1}^l\bE\left(f(\#\C_n(\x_i,\omega))\right)\,,\\
\end{eqnarray*}
where the first inequality follows from (\ref{ineqReimer}). Finally, we may let $n$ tend to infinity, and then use the Lebesgue's monotone convergence theorem for the right hand side and Fatou's lemma for the left hand side.    
\begin{flushright}\endproof\end{flushright}   

Let $\partial_\infty\A$ denote the lattice boundary of $\A$  with respect to the $l_{\infty}$-norm, and define 
$$\bar\A:=\A\cup\partial_\infty\A\,\mbox{ and }\,\C_\calY(\A):=\bar\A\cup\left\{\bigcup_{\C\in\calC_\calY(\A)}\bar\C\right\}\,.$$ 
\begin{lemma}\label{cluster}
If $\xi=\xi(\rho):=\bE e^{\#\C_\0}<\infty$ and $r\geq 2(\log \alpha+3^d+\log\xi)s$ then
$$\bP\left(\max_{\A\in\Phi_{\leq s}}\#\C_\calY(\A)\geq r\right)\leq e^{-r/2}\,.$$  
  \end{lemma}
  \proofof{Lemma \ref{cluster}} Notice that
 $$\#\C_\calY(\A)\,\leq\, \#\bar\A+\sum_{\C\in\calC_\calY(\A)}\#\C\,\leq\, 3^d\#\A+\sum_{\C\in\calC_\calY(\A)}\#\C\,.$$
By Markov's inequality and Lemma \ref{lemcorrelneg}, if $\A\in\Phi_{\leq s}$ then 
$$\bP\left(\sum_{\C\in\calC_\calY(\A)}\#\C\geq y\right)\leq e^{-y} \bE\left( e^{\sum_{\C\in\calC_\calY(\A)}\#\C}\right)\leq e^{-y}\left(\bE e^{\#\C_\0}\right)^s=e^{-y}\xi^s\,.$$
Hence
\begin{eqnarray}
\nonumber\bP\left(\max_{\A\in\Phi_{\leq s}}\#\C_\calY(\A)\geq r\right)&\leq&\bP\left(\max_{\A\in\Phi_{\leq s}}\sum_{\C\in\calC_\calY(\A)}\#\C\geq (r-3^ds)\right)\\
\nonumber&\leq&\alpha^s  e^{-(r-3^ds)}\xi^s\\
\nonumber&=&\exp\left\{ -r + s(\log\alpha+3^d+\log\xi)\right\}\\
\nonumber&\leq& e^{-r/2} \,,
\end{eqnarray}
whenever $r\geq 2(\log \alpha+3^d+\log\xi)s$.
\begin{flushright}\endproof\end{flushright}   

\subsection{Domination by product measures}  
Let $\calX=\{X_\z\,:\,\z\in\bZ^d\}$ be a collection of random variables that take values $0$ and $1$ and which satisfy the following conditions: (i) for each pair $\A,\B\in\bZ^d$ such that all sites in $\A$ are at distance greater than $k$ from all sites in $\B$ (in the sup-norm sense), the collections of random variables $\{X_{\z}\,\,:\,\,\z\in\A\}$ and $\{X_{\z}\,\,:\,\,\z\in\B\}$ are independent; (ii) $\inf_{\z\in\bZ^d}\bP\left(X_\z=1\right)\geq p$. In this case we say that $\calX$ is a $k$-dependent random field whose marginals are at least $p$, and denote $\calC(d,k,p)$ the class of all such fields. 

Let $\calY$ and $\calX$ be two random fields. We say that $\calY$ dominates $\calX$ from below if there is a coupling (joint realization) between $\calY$ and $\calX$ such that $Y_\z\leq X_\z$ for all $\z\in\bZ^d$. We refer to Liggett's book \cite{L85} for more details in stochastic domination and couplings. Theorem 0.0 of Liggett, Schonmann and Stacey \cite{LSS97} states that when $p$ is close to $1$, the random fields in $\calC(d,k,p)$ are dominated from below by an i.i.d. random field $\calY$ with density $\rho=\rho(d,k,p)$. Further, one can make $\rho$ arbitrarily close to $1$ by taking $p$ close enough to $1$. 

\begin{lemma}\label{resume}
Let $\calX\in\calC(d,k,p)$. There exists $\bar{p}\in(0,1)$ such that for all $p\in[\bar{p},1]$ if $s\geq 2r$ then 
$$\bP\left(\min_{\A\in\Phi_{\geq s}}\sum_{\z\in\A} X_\z \leq r\right)\leq e^{-s}\,.$$
 Further, there exists a constant $c_0>0$ such that if $r\geq c_0s$ then
 $$\bP\left(\max_{\A\in\Phi_{\leq s}}\#\C_\calX(\A)\geq r\right)\leq e^{-r/2}\,.$$  
\end{lemma}
\proofof{Lemma \ref{resume}}  We note that if $\calY$ dominates $\calX$ from below then
$$\sum_{\z\in\A} Y_\z \leq \sum_{\z\in\A} X_\z\,\,\mbox{ and } \,\,\#\C_\calX(\A)\leq \#\C_\calY(\A)\,,$$
for any lattice animal $\A$. Hence, Lemma \ref{resume} follows by combining Theorem 0.0 in \cite{LSS97} together with Lemma \ref{opendensity} and Lemma \ref{cluster}. 
\begin{flushright}\endproof\end{flushright}  

The proof of the first part of Lemma \ref{resume} (in the $k$-dependent set up) could be done directly without using \cite{LSS97}. One needs to notice that given any set of $s-r$ boxes, we can pick a subset of independent boxes of size $(s-r)/k^d$, and then use the same argument as before.    
However, the proof of the second part is more delicate and it is not clear (for the author) that it could be easily adapted to the $k$-dependent situation. 

\subsection{The block argument}
For each $\z\in\bZ^d$, $L>0$ and $s\in\{j/2\,:\,j\in\bN\}$ let 
$$B_\z^{s,L}:= L\z+[-sL,sL]^d\,.$$
Given a locally finite set $\calN\subseteq\bR^d$, we say that a (square) box $B$ is a $\calN$-\emph{full} box if cutting it
regularly into $(4\lceil\sqrt{d}\rceil+1)^d$ sub-boxes, each one of these boxes
contains at least one point of the set $\N$. Let  
$$B(\A):=\bigcup_{\z\in\A}B_{\z}^{1/2,L}\,\mbox{ and }\tilde{B}(\A)=\left\{\x\in\bR^d\,:\,\exists\,\,\y\in B(\A)\,\mbox{ s.t. }\,\|\x-\y\|_2\leq L/2\right\}\,.$$
\begin{lemma}\label{fullcontrol1}
Let $\A$ be a finite and connected subset of $\bZ^d$ and assume that $B_\z^{1/2,L}$ is a $\calN$-full box for all $\z\in\partial_\infty\A$.  If $C_\v\in\V$ and $C_\v\cap B(\A)\neq\emptyset$  then $C_\v\subseteq \tilde B(\A)$. 
\end{lemma}

\proofof{Lemma \ref{fullcontrol1}} Assume that $C_\v\cap B(\A)\neq\emptyset$ but $C_\v\not\subseteq \tilde B(\A)$. Then there will exist $\x_1,\x_2\in C_\v$ such that 
$$\|\x_1-\x_2\|_2\geq L/2\,.$$
On the other hand, $B_\z^{1/2,L}$ is a full box for all $\z\in\partial_\infty\A$. By picking $\x_1$ and $\x_2$ in the (euclidean) boundary of $B(\A)$ and $\tilde B(\A)$, respectively, this implies that there exist $\v_1,\v_2\in\calN$ such that  
$$\|\v_1-\x_1\|_2\leq \frac{\sqrt{d}}{4\lceil\sqrt{d}\rceil+1}L\,\,\mbox{ and }\,\,\|\v_2-\x_2\|_2\leq \frac{\sqrt{d}}{4\lceil\sqrt{d}\rceil+1}L\,$$ 
(the right hand side of the inequality is the length of the diagonal of a subsquare). However, 
$$\|\v-\x_1\|_2\leq\|\v_1-\x_1\|_2\,\,\mbox{ and }\,\,\|\v-\x_2\|_2\leq\|\v_2-\x_2\|_2\,,$$
and hence
\begin{eqnarray}
\nonumber\frac{1}{2}L&\leq& \|\x_1-\x_2\|_2\\
\nonumber&\leq& \|\x_1-\v\|_2+\|\x_2-\v\|_2\\
\nonumber&\leq&  \|\x_1-\v_1\|_2+\|\x_2-\v_2\|_2\leq \frac{2\sqrt{d}}{4\lceil\sqrt{d}\rceil+1}L\,,
\end{eqnarray}
which yields to a contradiction since $4\lceil\sqrt{d}\rceil+1>4\sqrt{d}$.
\begin{flushright}\endproof\end{flushright}   

\begin{lemma}\label{fullcontrol2}
Under \reff{comparaison},
$$\sup_{\z\in\bZ^d}\bP\left(\mbox{$B_\z^{1/2,L}$ is not a $\calN$-full box}\right)\leq (4\lceil\sqrt{d}\rceil+1)^d\exp\left\{-c_\mu^{-1}\left(\frac{L}{4\lceil\sqrt{d}\rceil+1}\right)^d\right\}\,.$$
\end{lemma}

\proofof{Lemma \ref{fullcontrol2}} Cut $B_\z^{1/2,L}$ regularly into $(4\lceil\sqrt{d}\rceil+1)^d$ sub-boxes, so that 
$$B_\z^{1/2,L}=\bigcup_{i=1}^{(4\lceil\sqrt{d}\rceil+1)^d} B_{\z,i}\,\,\,\mbox{ and }\,\,\,\bP\left(\#_\calN B_{\z,i}=0\right)\leq\exp\left\{-c_\mu^{-1}\left(\frac{L}{4\lceil\sqrt{d}\rceil+1}\right)^d\right\}\,.$$
Hence
\begin{eqnarray}
\nonumber\bP\left(\mbox{$B_\z^{1/2,L}$ is not a $\calN$-full box}\right)&\leq& \sum_{i=1}^{(4\lceil\sqrt{d}\rceil+1)^d}\bP\left(\#_\calN B_{\z,i}=0\right)\\
\nonumber&\leq& (4\lceil\sqrt{d}\rceil+1)^d \exp\left\{-c_\mu^{-1}\left(\frac{L}{4\lceil\sqrt{d}\rceil+1}\right)^d\right\}\,.
\end{eqnarray}
\begin{flushright}\endproof\end{flushright}   

\section{Proof of the Theorems}\label{Pr}
\paragraph{\bf Proof of \reff{e1t1}} For any connected set $C\subseteq\bR^d$, 
$$C\subseteq \bigcup_{\z\in\A(C)} B_\z\,,$$
and thus 
$$r\leq\#_\N\P\leq \sum_{\z\in\A(\P)} \#_\N B_\z=\sum_{\z\in\A(\P)} N_\z\,,$$
if $\P\in\Pi_{\geq r}$. Therefore, by Lemma \ref{Poisson}, 
\begin{equation}\label{rhs}
\bP\left( \min_{\P\in\Pi_{\geq r}}\#\A(\P)\leq s\right)\leq \bP\left(\max_{\A\in\Phi_{\leq s}}\sum_{\z\in\A} N_\z\geq r\right)\leq e^{-r/2} \,,
\end{equation}
whenever $r\geq b_1s$ and $b_1:= 2(\log \alpha +c_\mu(e-1))$. 

\paragraph{\bf Proof of \reff{e2t1}} For each $L\geq 1$ let
$$\A_L(\P):=\left\{\z\in\bZ^d\,\,:\,\,B_\z^{1/2,L}\cap \P\neq\emptyset\,\right\}\,,$$
and recall that $\A(\P)=\A_1(C)$. Then, for any $L\geq 1$,
\begin{equation}\label{L}
\#\A_L(\P)\leq \#\A_1(\P)\leq L^d \#\A_L(\P)\,.
\end{equation}
Consider the  non-homogeneous $3$-dependent percolation scheme $\calX^L$ in $\bZ^d$ defined by
$$X_\z^L:=\mathbf{1}\left\{B_{\z'}^{1/2,L}\,\,\mbox{ is a full box }\forall\,\,\z'\mbox{ s.t. }\,\,\|\z'-\z\|_{\infty}\leq 1\right\}\,.$$
($\mathbf{1}$ denotes the indicator function of an event.) If $X_\z^L=1$ we say that $\B_{\z}^{1/2,L}$ is a good box. By Lemma \ref{fullcontrol2}, $\calX^L\in\calC(d,3,p_L)$ where
$$1-p_L:=\sup_{\z\in\bZ^d}\bP\left(X^L_\z=0\right)\leq 3^d(4\lceil\sqrt{d}\rceil+1)^d\exp\left\{-c^{-1}_\mu\left(\frac{L}{4\lceil\sqrt{d}\rceil+1}\right)^d\right\}\,,$$
By Lemma \ref{resume}, if we pick $L_0$ such that
$$3^d(4\lceil\sqrt{d}\rceil+1)^d\exp\left\{-c^{-1}_\mu\left(\frac{L_0}{4\lceil\sqrt{d}\rceil+1}\right)^d\right\}\leq1- \bar{p}\,$$
then $p_{L_0}\in [\bar{p},1]$ and 
\begin{equation}\label{super}
\bP\left(\min_{\A\in\Phi_{\geq s}}\sum_{\z\in\A} X^{L_0}_\z \leq r\right)\leq e^{-s}\,.
\end{equation}
whenever $s\geq 2r$. Now, let
$$\calS_{\calX^L}(\P):=\{\z\in\A_L(\P)\,:\,X_\z^{L}=1\}\,.$$  
Notice that there exists at least one set $\calS'_\calX\subseteq\calS_\calX$ such that $|\z-\z'|_{\infty}\geq 2$ for all $\z,\z'\in\calS'_\calX$ and $k=\#\calS'_\calX\geq \#\calS_\calX /3^d$. Now, write $\calS'_\calX=\{\z_1,\dots,\z_k\}$. By Lemma \ref{fullcontrol1}, if $\z_i,\z_j\in\calS'$ and $C_{\v_i}\cap B_{\z_i}^{1/2,L_0}\neq\emptyset$ and $C_{\v_j}\cap B_{\z_j}^{1/2,L_0}\neq\emptyset$ then $\v_i\neq\v_j$, and thus
$$\#_\N\P\geq k \geq \frac{\#\calS_\calX}{3^d}\geq\frac{\sum_{\z\in\A(\P)}X_\z^L}{3^d}\,,$$ 
By \reff{L}, this shows that 
\begin{eqnarray}
\nonumber\bP\left(\max_{\P\in\Pi_{\leq r}}\#\A(\P)\geq s\right)&\leq&\bP\left(\max_{\P\in\Pi_{\leq r}}\#\A_{L_0}(\P)\geq s/L_0\right)\\
\nonumber&\leq& \bP\left(\min_{\A\in\Phi_{\geq s/L_0}}\sum_{\z\in\A} X^{L_0}_\z \leq 3^d \,r\right)\leq e^{-s/L_0}\,,
\end{eqnarray}
whenever $s\geq 2(L_0 3^d)r$.
   
\paragraph{\bf Proof of \reff{e2t2}} By Lemma \ref{fullcontrol1} (recall the definition of $\C_{.}$ in Lemma \ref{cluster}), 
$$\P(\A)\subseteq B\left(\C_{\calX^L}(\A)\right)\,$$
and
$$\#_\N\P(\A)\,\,\leq\,\, \#_\N B\left(\C_{\calX^L}(\A)\right)=\sum_{\z\in\C_{\calX^L}(\A)}N_\z\,.$$
which yields to
$$\bP\left(\max_{\A\in\Phi_{\leq s}}\#_\N\P(\A)\geq r\right)\leq \bP\left(\max_{\A\in\Phi_{\leq cr}}\sum_{\z\in\A} N_\z\geq r\right)+ \bP\left(\max_{\A\in\Phi_{\leq s}}\#\C_{\calX^L}(\A)\geq cr\right)\,,$$
for any $c>0$. By Lemma \ref{Poisson} and Lemma \ref{cluster},
$$\bP\left(\max_{\A\in\Phi^{cr}}\sum_{\z\in\A} N_\z\geq r\right)\leq e^{- r/2}\,\,\,\mbox{ and }\,\,\,\bP\left(\max_{\A\in\Phi_{\leq s}}\#\C_{\calX^L}(\A)\geq cr\right)\leq e^{-cr/2}\,,$$
for $c=[2(\log\alpha+c_\mu(e-1))]^{-1}$ and $r\geq (c_0/c) s$ (notice that $r\geq (\log\alpha+c_\mu(e-1))cr=r/2$), which shows that 
$$\bP\left(\max_{\A\in\Phi_{\leq s}}\#_\N\P(\A)\geq r\right)\leq e^{- r/2}+ e^{-cr/2}\,,$$
and finishes the proof of \reff{e2t2}. 

\section{A two dimensional modified Poisson model}\label{Mo}
In \cite{P07} a random set is constructed from a realization of a two-dimensional homogeneous  Poisson random set $\N$ as follows. Order the points of $\bZ^{2}$ in some arbitrary fashion, say $\bZ^2:=\{\u_1,\u_2, \dots\}$. Fix $\delta\in(0,1)$ and $n\geq 1$. For each $k\geq1$ let 
\[
B_k^{n}:=B_{\u_k}^{1/2,n^{\delta}}\,. 
\]
Divide $B_k^n$ into $36$ sub-boxes (as before) of the same length $n^\delta/6$, say  $B_{k,1}^n\dots B_{k,36}^n$. Now we construct the modified random set $\calN(n):=\calN(n,\calN)$ (whose distribution will also depend on $n$ and $\delta$) by changing the original Poisson random set $\calN$ inside each $B_{k,j}^n$, as follows. 
\begin{enumerate}
\item If $1\leq \#_{\calN}B_{k,j}^{n}\leq n^{2\delta}$ then set $B_{k,j}^{n}\cap\calN(n):=B_{k,j}^{n}\cap\calN$;
\item If $\#_{\calN}B_{k,j}^n> n^{2\delta}$ then set $\calN(n)$ by uniformly selecting  $n^{2\delta}$ points from $B_{k,j}^{n}\cap\calN$.
\item If $\#_{\calN}B_{k,j}^{n}=0$ then set $B_k^{n}\cap\calN(n)$ by adding an extra point  uniformly distributed on $B_{k,j}^{n}$.
\end{enumerate} 
In a few words, we tile the plane into boxes $B_k^n$ of size $n^\delta$ and we insist that each tile is a full box, and that no tile contains more than $36 n^{2\delta}$ Poisson points. We make the convention $\calN(\infty)=\calN$ and denote by $\D(n)$ the Delaunay Triangulation based on $\calN(n)$. It is clear that method can be applied in this set up for each fixed $n$. However, we want to emphasize that it allows us to do so simultaneously for all $n$  sufficiently large. 
\begin{theorem}\label{t3}
Theorems \ref{t1} and \ref{t2}, as well as Corollaries \ref{c1} and \ref{c2}, hold in the modified model $\N(n)$, where the constants $b_j$ does not depend $n\geq 1$.   
\end{theorem}
\proofof{Theorem \ref{t3} (outline)} The proof of \reff{e1t1}, in the $\N(n)$ context, uses that 
$$\#_{\N(n)}B_\z\leq \#_\N B_\z+2^d\,,$$
(the small box could intersect at most $2^d$ boxes $B_{k,j}^n$) and hence (see \reff{rhs})
$$\bP\left( \min_{\P\in\Pi_{\geq r}}\#\A(\P)\leq s\right)\leq \bP\left(\max_{\A\in\Phi_{\leq s}}\sum_{\z\in\A} N_\z\geq r-2^ds\right)\,,$$     
which implies \reff{e1t1} (by Lemma \ref{Poisson}). 

The proof of \reff{e2t1} and \reff{e2t2} were essentially based on the fact that we can compare the original problem to a site percolation scheme with minimal marginal density $p_L\to1$ as $L\to\infty$. The same comparison method works here as soon as we prove that 
$$\inf_{n}p_L(n)\to 1\,\,\mbox{ as }\,\,L\to \infty\,,$$
or, equivalently, 
$$\lim_{L\to\infty}\sup_{n\geq 1}\bP\left(\mbox{$B_\z^{1/2,L}$ is not a $\calN(n)$-full box}\right)=0\,.$$
Notice that, if $L\geq 2n^\delta$ then a box $B$ of size $L/6$ must contain at least one box  $B_{k,j}^{n}$ of length $n^\delta/6$ (for some $k$ and $j$) which means that $B$ is a $\N(n)$-full box. 
On the other hand, if $L\leq 2n^{\delta}$ then the probability of the event that $B_\z^{1/2,L}$ is not a $\N(n)$-full box is bounded by the probability of the event that $B_\z^{1/2,L}$ is not a $\N$-full box  plus the probability of the event that $\N\cap B_\z^{1/2,L}\neq\N(n)\cap B_\z^{1/2,L}$, which decays to $0$ as $\max\{e^{-L^2},e^{-n^{2\delta}}\}\leq e^{-L^2}$. The rest of the proof follows mutatis-mutandis the method applied to prove \reff{e2t1} and \reff{e2t2}.  
 \begin{flushright}\endproof\end{flushright}  
 
The techniques developed in this paper also fit to study the minimal density of open edges among all self-avoiding paths $\gamma\in\Gamma_{\geq r}$, in the two dimensional bond percolation model in the Delaunay triangulation. This model is constructed by attaching i.i.d. Bernoulli random variables $\tau_\e$ with parameter $p\in[0,1]$ to edges $\e\in\D$. Let $p_c^*\in (0,1)$ be the critical probability for the bond percolation model in the $\N$-Voronoi tessellation \cite{P06}. We denote by $(\calN(n),\tau)$ the the bond percolation model associated to the modified random set $\N(n)$. 
\begin{theorem}\label{thm:reward}
If $\bP(\tau_\e=0)<1-p_c^*$ then there exists $n_0\geq 0$ and $b_9,b_{10}>0$ such that, for all $n\geq n_0$, if $\,r\geq b_9 s\,$  then
\begin{equation}\label{reward}
\bP\left(\min_{\gamma\in\Gamma_{\geq r}}\sum_{\e\in\gamma}\tau_\e\leq s\right)\leq 3 e^{-b_{10} r}\,.
\end{equation}
\end{theorem}

\proofof{Theorem \ref{thm:reward} (outline)}  Let $\Gamma_\z^L$ be the collection of all self-avoiding paths $\gamma=(\v_0,\v_1,\dots,\v_l)$ in $\D$ such that $C_{\v_0}\cap\partial B_\z^{1/2,L}\neq\emptyset$, $C_{\v_l}\cap\partial B_\z^{3/2,L}\neq\emptyset$ and $\v_j\in B_\z^{3/2,L}\setminus B_\z^{1/2,L}$ for $j=1,\dots,l-1$. We say that $B_\z^{1/2,L}$ is a $(\N,\tau)$-good box if the following holds: 
\begin{enumerate}
\item $B_{\z'}^{1/2,L}$ is a $\N$-full box for all $\z'\in\bZ^2$ s.t. $|\z'-\z|_{\infty}\leq 2$;
\item $\sum_{\e\in\gamma}\tau_\e\geq 1$ for all $\gamma\in\Gamma_\z^L$ . 
\end{enumerate}
We have seen that the probability of the event (1) is uniformly bounded for all $n\geq 1$. For $n=\infty$ (we are in classical Poisson model), Lemma 1 in \cite{P06} implies that the probability of the event (2) goes to $0$ as $L\to\infty$, if $\bP(\tau_\e=0)<1-p_c^*$. For a finite $n\geq 1$ the probability that $\N\neq\N(n)$ inside $B_\z^{1/2,L}$ is of order $L^2e^{-2n^\delta}$. Thus,
\begin{eqnarray}
\nonumber\bP\left(\mbox{(2)}\mid \mbox{(1)}\right)&\leq& \bP\left(\mbox{(2) for $n=\infty$}\mid \mbox{(1)}\right)
+cL^2e^{-2c'n^\delta}\\
\nonumber&\leq&  \bP\left(\mbox{(2) for $n=\infty$}\mid \mbox{(1)}\right)+cL^2e^{-2c'L^\delta} \,,
\end{eqnarray}
for $L\leq n$ ($c$ and $c'$ are constants). This means that, given any $p\in(0,1)$ there exists $L_0>0$ such that for all $n,L\geq L_0$,
\begin{equation}\label{marginal}
 \bP\left( B_\z^{1/2,L}\,\,\mbox{ is a $(\N(n),\tau)$-good box }\right)\geq p\,.
 \end{equation}
Consider the homogeneous percolation scheme $\calZ^L$ defined by
$$Z_\z^{L}:=\mathbf{1}\left\{B_\z^{1/2,L}\,\,\mbox{ is a $(\N(n),\tau)$-good box }\,\,\right\} \,.$$
By Lemma 3 in \cite{P06}, it is a $5$-dependent percolation scheme. Together with \reff{marginal}, this implies that $\calZ^{L_0}\in\calC(d,5,p)$ for all $n\geq L_0$, if $\bP\left(\tau_\e=0\right)< 1-p_c^*$. Now, let
$$\calS_{\calZ^L}(\gamma):=\{\z\in\A_L(\gamma)\,:\,Z_\z^{L}=1\}\,.$$  
Notice that there exists at least one set $\calS'\subseteq\calS$ such that $|\z-\z'|_{\infty}\geq 4$ for all $\z,\z'\in\calS'$ and $k=|\calS'|\geq |\calS|/4^d$. Now, write $\calS'=\{\z_1,\dots,\z_k\}$. By Lemma  \ref{fullcontrol1}, one can find disjoint pieces of $\gamma$, say $\gamma_1,\dots,\gamma_k$, such that, for $i=1,\dots,k$, $\sum_{\e\in\gamma_i}\tau_\e\geq 1$. Hence,
$$
\sum_{\e\in\gamma}\tau_\e\geq\sum_{i=1}^{k}\left( \sum_{\e\in\gamma_i}\tau_\e\right)\geq k=|\calS'|\geq \frac{|\calS|}{4^d}=\frac{\sum_{\z\in\A(\gamma)}Z_\z^{L}}{4^d}\,,
$$
which shows that 
\begin{eqnarray}
\nonumber\bP\left(\min_{\gamma\in\Gamma_{\geq r}}\sum_{\e\in\gamma}\tau_\e\leq s\right)&\leq& \bP\left(\min_{\gamma\in\Gamma_{\geq r}} |\A_{L_0}(\gamma)|\leq \frac{r}{b_1}\right)+\bP\left(\min_{\A\in\Phi_{\geq  r/b_1}}\sum_{\z\in\A} Z_\z^{L_0}\leq 4^d s\right)\\
\nonumber &\leq& e^{-r/2}+ 2 e^{- r/b_1}\,,
\end{eqnarray}
whenever $r\geq 2b_1(4^d)s$. In the last line, we have used \reff{e1t1} (in the modified context) and Lemma \ref{resume} (by choosing $L_0$ large enough).
\begin{flushright}\endproof\end{flushright}

\end{document}